\documentclass[]{article}
\usepackage{amsfonts}
\usepackage{theorem}

\newtheorem{theorem}{Theorem}
\newtheorem{corollary}{Corollary}

\newcommand{\kocka}{$\begin{array}{c}
\hspace*{-0.55em}\sqcap \hspace*{-0.60em}\\[-0.4em] \hline
\multicolumn{1}{c}{\hspace*{-0.60em}}\\[-0.8em]
\end{array}
$}

\usepackage{enumerate}

\begin{document}

\centerline{\bf On a problem of A. Nagy concerning permutable semigroups} 
\centerline{{\bf satisfying a non-trivial permutation identity}\footnote{Research supported by the Hungarian NFSR grant No T042481}}

\bigskip
\centerline{Attila De\'ak}

\bigskip

\centerline{\it Communicated by M. B. Szendrei}

\bigskip

\begin{abstract}
We say that a semigroup $S$ is a permutable semigroup if, for every
congruences $\alpha$ and $\beta$ of $S$, $\alpha \circ \beta =
\beta \circ \alpha$.
In [4], A. Nagy showed that every
permutable semigroup satisfying an arbitrary non-trivial permutation identity
is medial or an ideal extension of a rectangular band by a non-trivial commutative
nil semigroup. The author raised the following problem: Is every permutable semigroup
satisfying a non-trivial permutation identity medial?  In present paper we give a positive
answer for this problem.
\end{abstract}

\bigskip

A semigroup $S$ is called a {\it permutable semigroup} if
the congruences of $S$ commute with each other, that is,
$\alpha \circ \beta =\beta \circ \alpha$ is satisfied
for every congruences $\alpha$ and $\beta$ of $S$.

\begin{theorem}\label{tet1}
Let $S$ be an ideal extension of a rectangular band $K$ by a non trivial commutative nil semigroup $N$.\\
If $S$ is permutable and satisfies a non-trivial permutation identity then either $S$ is medial or $N$ is nilpotent.
\end{theorem}

\noindent
{\bf Proof}. Assume that a permutable semigroup $S$ satisfies a permutation identity
$$x_1x_2\ldots x_n=x_{\sigma(1)}x_{\sigma(2)}\ldots x_{\sigma(n)},\ (\sigma \neq id).$$
We can suppose that $n>2$.
Let $l$ and $k$ be non-negative integers such that\\
$$\sigma(1)=1,\sigma(2)=2,\cdots,\sigma(l)=l,\sigma(l+1)\neq l+1$$
and
$$\sigma(n-k)\neq n-k,\sigma(n-k+1)=n-k+1,\cdots,\sigma(n)=n.$$
If $l=k=0$ then, by Theorem 1 of [4], $S$ is commutative.
Consider the case when $l>0$ or $k>0$. Let $j$ be $\max\{l,k\}$.

Assume that $j=1$. Then let $\beta$ be a relation on $S$ such that $a\,\beta\,b\Leftrightarrow uabv=ubav$ for every $u,v\in S$. It is clear that $\beta$ is reflexive and symmetric. We show that $\beta$ is transitive. Let $u, v, a, b, c\in S$ be arbitrary elements such that $a\,\beta\,b$ and $b\,\beta\,c$. If $ac\in S-K$ then $ac=ca$, because $N$ is commutative. Thus $a\,\beta\,c$.
Assume that $ac\in K$. Then $ca\in K$ and $$ac=ac(xac)ac=acxac,\ ca=ca(xca)ca=caxca$$
for every $x\in S$. Thus $$uacv=(ua)cb(acv)=uab(cacv)=ubacacv=ubacv=uabcv=uacbv.$$ Similarly,
$$ucav=ubcav=ucbav=ucabv.$$ From these results it follows that $$uacv=ub^pab^qcb^rv$$ for every non-negative integers $p, q, r$. Similarly,
$$ucav=ub^pcb^qab^rv$$ for every non-negative integers $p, q, r$.
Let $\hat{l},\ \hat{k}$ be such that $\sigma(l+1)=\hat{l},\ \sigma(n-k)=\hat{k}$. It is obvious that $\hat{l}\geq l+2$ and $ \hat{k} \leq n-k-1$.
So $$uacv=uacb^{n-2}v=ub^lab^{\hat{l}-l-2}cb^{n-\hat{l}}v=ub^lcb^tab^sv=ucav,$$ where $s, t$ are non-negative integers such that $s+t+l=n-2$.
Hence $\beta$ is an equivalence on $S$. Let $x,y\in K$. We will show that $x\,\beta\,y$.
Let $u,v\in S$ be arbitrary elements. Then
$$uxyv=ux^{\hat{l}-l-1}y^{n-\hat{l}}yv=uyxyv=uyv=uy^{\hat{k}}xy^{n-\hat{k}-k-1}v=uy^{n-k-1}xv=uyxv.$$ Since, for every $a\in N$, there is a positive integer $m$ such that $a^m\in K,\ a\,\beta\,a^m$, we get that $S$ is a $\beta$ class. Hence $S$ is medial.

Consider the case when $j>1$.
Let $\alpha$ be an equivalence of $S$ such that 
$a\,\alpha\,b\Leftrightarrow \forall x\in K\ xa=xb$ and $ax=bx$.\\
It is obvious that $x\,\alpha\,y$ implies $x=y$ for every $x,y\in K$.
We show that $\alpha$ is also a congruence.
Let $y,a,b\in S$ be arbitrary elements. Assume $a\,\alpha\,b$ then $$x(ay)=(xa)y=(xb)y=x(by)$$ and $$(ay)x=a(yx)=b(yx)=(by)x$$ and so $ay\,\alpha\,by$.
We can prove $ya\,\alpha\,yb$ in a similar way.
Assume that $N$ is not nilpotent. Let $a_1,a_2,\cdots,a_j\in S-K$ be such that $a_1a_2\ldots a_j\in S-K$.
We prove that $a_1a_2\ldots a_j\,\alpha\,a_1a_2\ldots a_j^m$ where $m$ is the least positive integer such that $a_j^m\in K$.
Since $a_1a_2\ldots a_j\in S-K$ then 
$ a_{i_1}a_{i_2}=a_{i_2}a_{i_1}$ for every $i_1,i_2\in \{1,2,\cdots ,j\}$.
Let $l^{'}$ and $k^{'}$ be such that $\sigma(l+1)=l^{'}$ and $\sigma(n-k)=k^{'}$.
Then, for every $x\in K$, $$a_1a_2\ldots a_jx=a_1a_2\ldots a_jx^{l^{'}-l-1}a_jx^{n-l^{'}}x=a_1a_2\ldots a_ja_jx^{n-l-1}x=a_1a_2\ldots a_j^2x$$
and so $$a_1a_2\ldots a_jx=a_1a_2\ldots a_j^mx.$$
Using also $a_{i_1}a_{i_2}=a_{i_2}a_{i_1}\ (i_1,i_2\in \{1,2,\cdots ,j\})$,
we have $$xa_1a_2\ldots a_j=xa_1a_2\ldots a_j^m$$ in a similar way.
From these results we have $$a_1a_2\ldots a_j\,\alpha\,a_1a_2\ldots a_j^m.$$
From (2) of Lemma 3 of [3], we get $\mid K \mid =1$, and so $S$ is commutative.
\hfill\kocka

\begin{theorem}\label{tet2}
If $S$ is an ideal extension of a rectangular band $K$ by a non trivial commutative nilpotent semigroup $N$, and $S$ is permutable then S is commutative.
\end{theorem}

\noindent
{\bf Proof}.
Assume that a permutable semigroup $S$ is
an ideal extension of a rectangular band $K$ by a non trivial commutative nilpotent semigroup $N$.
By Lemma 2 of [4], $N$ is a $\Delta$-semigroup. Then,
by Lemma 2 of [5], $N$ is finite ciclyc.
Let $N=\{a,a^2,\cdots,a^n\}$ where $a^n\in K$. Assume $a^n=(1,1)$ and $K=L\times R$, where $L$ is a leftzero,
$R$ is a rightzero semigroup.
For every $(x_1,x_2) \in K$, $a(x_1,x_2)=(\phi_a(x_1,x_2),x_2)$
and $(x_1,x_2)a=(x_1,\psi_a(x_1,x_2))$.\\
If $(x_1,x_3)\in K$, then $$a(x_1,x_3)=a(x_1,x_2)(x_2,x_3)=(\phi_a(x_1,x_2),x_2)(x_2,x_3)=(\phi_a(x_1,x_2),x_3)$$ and so $\phi_a(x_1,x_2)=\phi_a(x_1,x_3)$. So $\phi_a$ is independent from the second variable. Similarly, $\psi_a$ is independent from the first variable. From these results we get that, for every $(x_1,x_2)\in K$, $a(x_1,x_2)=(\phi_a(x_1),x_2)$ and $ (x_1,x_2)a=(x_1,\psi_a(x_2))$.
We will prove that $\phi_a(x)=1$ for every $x\in L$. It is obvious that $\underbrace{\phi_a\circ\ldots\circ\phi_a}_{\mbox{n}}(x)=\phi_a^n(x)=1$ for every $x\in L$.
Let $j$ be such that $\phi_a^j(x)=\phi_a^{j+1}(x)=\ldots =\phi_a^n(x)=1\ \forall x\in L$.
Assume that $\exists x\in L$ such that $\phi_a^{j-1}(x)\neq 1$. Let $$P:=\{x\in L:\,\phi_a^{j-1}(x)\neq 1\}$$
and
$$Q:=\{(x_1,x_2)\in K: x_1\in P\}.$$
It is clear that, for every $(x_1,x_2)\in K$, $a^r(x_1,x_2)\notin Q$,
for all positive integers $r$, because
$a^r(x_1,x_2)=(\phi_a^r(x_1),x_2)$ and $\phi_a^{j-1}(\phi_a^r(x_1))=\phi_a^{j+r-1}(x_1)=1$. 
Since $\phi_a(1)=1$ then $Q\neq K$.\\
Let $\alpha$ be an equivalence on $S$ such that, for arbitrary $x,y\in S$, $x\,\alpha\,y\Leftrightarrow x,y\in Q$ or $x,y\in S-Q$.
We show that $\alpha$ is a congruence of $S$.
Let $x,y,z\in S$ be arbitrary elements with $x\,\alpha\,y$.
Assume that $x,y\in K$. It is clear that if $z\in K$ then $xz\,\alpha\,yz$ and $zx\,\alpha\,zy$. Assume
$z=a^l\in S-K$. It is obvious that $xa^l\,\alpha\,ya^l$. Since $a^lx,a^ly\notin Q$ then $a^lx\,\alpha\,a^ly$ is obvious.
Consider the case when $x=a^m\in S-K$ and $y\in K$ where $m$ is a positive integer less than $n$.
For every $z\in K\ zx\,\alpha\,zy$, is obvious. Since $xz,y\notin Q$, we get that $xz\,\alpha\,yz$.
Consider the case when $x,y\in S-K$, that is $x=a^p$ and $y=a^q$, where $p,q$ are positive integers less than $n$.
For all $z\in S$, $zx\,\alpha\,zy$ and $xz\,\alpha\,yz$ are obvious.
So $\alpha$ is a congruence, indeed.
There are two $\alpha$ classes on $S$: $Q$ and $S-Q$.
It is clear that $a^n\in K-Q$ and so $Q\neq K$.
From (2) of Lemma 3 of [3], $K$ have to be union of $\alpha$ classes or
subset of an $\alpha$ class.
Since $(S-Q)\cap K\neq 0$, $(S-Q)\cap (S-K)\neq 0$ and $Q\subset K$,
then $Q$ have to be an empty set! Since $\mid Q\mid \geq \mid P\mid$, then $0=\mid Q\mid\geq\mid P\mid\geq 0$.
That means $\mid P\mid =0$, that is $\phi_a(L)=\{1\}$.
We can prove $\psi_a(R)=\{1\}$ in a similar way.
Consequently, for every $x\in K$, $ax=a^2x=\ldots =a^nx$ and $xa=xa^2=\ldots =xa^n$.\\
Let $\beta$ be an equivalence of $S$ such that $x\,\beta\,y\Leftrightarrow x=y$ or $x=a^l,\ y=a^k$ where $l,k$
are positive integers. Let $x,y,z\in S$ be arbitrary elements with $x\,\beta\,y$. If $x=y$ then $xz\,\beta\,yz$
and $zx\,\beta\,zy$ are obivous.
If $x\neq y$ then, with $x\,\beta\,y$, we get that $x=a^l$ and $y=a^k$, for some positive integers $l,k$. If $z=a^m$ then
$a^ma^l\,\alpha\,a^ma^k$ and $a^la^m\,\alpha\,a^ka^m$ are obvious. Let $x\in K$, then $xa^l=xa^n=xa^k$ and
$ a^lx=a^nx=a^kx$. Hence $\beta$ is a congruence.
It is clear that, for all $x,y\in K$, $x\,\beta\, y\Leftrightarrow x=y$ and, for all $x,y\in
\{a,a^2,\cdots ,a^n\}$, $x\,\beta\,y$.
So the $\beta$ classes of $S$ are:
$\{ a,a^2,\cdots ,a^n\}=N$ and for every $x\in S-N$, $x$ is a $\beta$ class.
From (2) of Lemma 3 of [3] $K$ have to be union of $\beta$
classes, or subset of a $\beta$ class. Since $\{ a,a^2,\cdots ,a^n\}\cap K=a^n$ and
$S=\{ a,a^2,\cdots ,a^{n-1}\}\cup K$, then $K={a^n}$. That is $\mid K\mid =1$.
And so $S$ is commutative.
\hfill\kocka

\begin{corollary}
Every permutable semigroup satisfying a non-trivial permutation identity is medial.
\end{corollary}

\noindent
{\bf Proof}. From Theorem 2 of [4] and Theorem 1 and Theorem 2 of this paper it is obvious.
\hfill\kocka

\bigskip

\noindent
Attila De\'ak, Department of Algebra, Institute of Mathematics,
Budapest University of Technology and Economics; e-mail: ignotus@math.bme.hu

\end{document}